# Analysis of a Leslie-Gower predator-prey model with Allee effect and fear effect

LIN CHEN[1] AND CHANGRONG ZHU[*]

**ABSTRACT:** In this paper, the dynamical behaviors of a Leslie–Gower predator–prey model with Allee effect and fear effect are studied. First, we use Blow-Up method to explore the stability of the original point. Analyze the local stability of the non-negative equilibriums of the system. Taking the growth rate of predators as the Hopf bifurcation parameters, and give the condition of the model Hopf bifurcation. Using the method of calculating the first Lyapunov coefficient to discuss the direction of the Hopf bifurcation near the weak center.

## 0 Introduction

In ecosystems, fear effects are prevalent between Predators and preys, where prey becomes fearful of the presence of predators and to some extent affects the population behavior of prey. In 2016, Wang et al[1]. found that prey exhibit anti-predation phenomena, i.e., changes in foraging, reproduction, and vigilance of prey populations, when they perceive the risk of being predated, based on a study of a terrestrial vertebrate. Based on this phenomenon, Wang et al. proposed the following predator-prey model with a fear effect, considering that the fear effect can negatively affect the birth rate of prey populations:

$$\begin{cases} \dot{u} = \dfrac{ru}{1+ku} - du - au^2 - g(u)v \\ \dot{v} = cg(u)v - mv \end{cases} \quad (0.1)$$

where $u$ and $v$ denote the population densities of food bait and predators at moment $t$, respectively, $r$ is the birth rate of predators, $d$ is the natural mortality preys, $g(u)$ denotes the functional response function, $c$ is the conversion rate of predators, and $m$ is the mortality rate of predators. $\dfrac{1}{1+ku}$ is the fear effect and $k$ is the fear intensity. Wang et al. performed a detailed kinetic analysis of the model and concluded that low levels of fear intensity can produce limit rings through Hopf bifurcation. After that, numerous scholars have analyzed and studied the effect of fear effect on system production, Wang J, Cai Y L, Fu S M et al[2] studied the effect of fear factor on predator-prey model with shelter. S. Pal et al[3] considered fear effects in models where predators have cooperative behavior and concluded that cooperation among predators exacerbates fear in prey populations. The rest of scholars[4,5] have also analyzed the dynamics of models with fear effects.

Allee effect is also one of the research focuses affecting ecosystems. For populations, too sparse a population density will affect the birth rate of the population, and too high a single population density will likewise increase interspecific competition, thus leading to a decline in population size, such a phenomenon is called Allee effect. In recent years, the Allee effect has been extensively studied by scholars due to the harshness of the natural environment and for the conservation of endangered animals. Arancibia-Ibarra[6] studied the dynamic shape of a model with Holling-Tanner functional function under the Allee effect, and Olivares[7] studied a class of predators with Allee effect was analyzed in the Leslie-Gower model.

Based on previous research results, this paper proposes a class of Leslie-Gower predator-prey systems with fear and Allee effects:

$$\begin{cases} \dfrac{dx}{dt} = \dfrac{rx}{1+\lambda y}\left(1-\dfrac{x}{K}\right)(x-m) - axy \\ \dfrac{dy}{dt} = sy\left(1-\dfrac{y}{hx}\right) \end{cases}$$

where $x$ and $y$ denote the population densities of prey and predators at time $t$, respectively, $x-m$ denotes the Allee effect, $\dfrac{1}{1+\lambda y}$ is the fear effect, and $\lambda$ is the fear intensity. $K$ is the carrying capacity, and $a$ is the intraspecific competition intensity factor of the bait population.

## 1 Model building

According to the predator-prey model proposed by Leslie and Gower, considering the Allee effect that prey populations have affects the population size; and that the fear effect of prey on predators decreases the birth rate of prey. Therefore, in conjunction with the Leslie-Gower model, model (1.1) is given as:

$$\begin{cases} \dfrac{dx}{dt} = \dfrac{rx}{1+\lambda y}\left(1-\dfrac{x}{K}\right)(x-m) - axy \\ \dfrac{dy}{dt} = sy\left(1-\dfrac{y}{hx}\right) \end{cases} \quad (1.1)$$

All parameters in model (1.1) are positive by the following scale transformation:

$$\bar{x} = \dfrac{1}{K}x, \quad \bar{y} = \dfrac{1}{hK}y, \quad \tau = rKt, \quad \bar{\lambda} = \lambda hK, \quad \bar{a} = \dfrac{ha}{r}, \quad \bar{s} = \dfrac{s}{rK} \circ$$

Let $\bar{x}$, $\bar{y}$, $\tau$, $\bar{\lambda}$, $\bar{a}$, $\bar{s}$ be $x$, $y$, $t$, $\lambda$, $a$, $s$. The model (1.1) can be transformed into (1.2):

$$\begin{cases} \dfrac{dx}{dt} = \dfrac{x}{1+\lambda y}(1-x)(x-m) - axy \\ \dfrac{dy}{dt} = sy\left(1-\dfrac{y}{x}\right) \end{cases} \quad (1.2)$$

The threshold of Allee effect is $m$ after dimensionless, and $0 < m < 1$ for strong Allee effect。

## 2 Dynamical form of the model equilibrium point

In this section, the dynamical form of the system will be studied. In the system (1.2), , X, Y denote the densities of the prey and predator populations, respectively, and are not defined on $x = 0$. From (1.2), the positive half-axis of the $x$ axis is the invariant set of the system. Considering the biological significance, the morphology of the system is discussed in the first quadrant and the positive half-axis of the $x$ axis in this paper.

### 2.1 Stability of the origin

By the assumptions $x(0) > 0$ and $y(0) \geq 0$. Since the system (1.2) is not meaningful at the origin and the dynamical properties of the origin cannot be analyzed by the Jacobi matrix, the blow-up method[8] is considered to discuss the stability of the origin.

**Theorem 1:** The origin of the system (1.2) is an unstable equilibrium point.

**Proof:** Using the blow-up method, let $x = u$ and $y = v$. Then the system (1.2) becomes the following system:

$$\begin{cases} \dfrac{du}{dt} = u\left(\dfrac{1}{1+\lambda uv}(1-u)(u-m) - auv\right) \\ \dfrac{dv}{dt} = v\left(s(1-v) - \dfrac{1}{1+\lambda uv}(1-u)(u-m) + auv\right) \end{cases} \quad (2.1)$$

Let $u = 0$, then from equation (2.1) we have that $v = 0$ or $v = \dfrac{s}{s-m}$. That is, two primary singularities $P_1 = (0,0)$ and $P_2 = \left(0, \dfrac{s}{s-m}\right)$ are obtained after blowing up the origin of the system (1.2). The Jacobi matrix of the system (2.1) at $P_1 = (0,0)$ is:

$$J(P_1) = \begin{pmatrix} -m & 0 \\ 0 & s+m \end{pmatrix}$$

The Jacobi matrix at $P_2 = \left(0, \dfrac{s}{s-m}\right)$ is given by:

$$J(P_2) = \begin{pmatrix} m & 0 \\ -\left(m+1+\dfrac{sm}{s-m}\right)\dfrac{sm}{s-m} + \dfrac{as}{s-m} & -\dfrac{s^2+m^2}{s-m} \end{pmatrix}$$

From the Jacobi matrices $J(P_1)$ and $J(P_2)$, when $s < m$, the system (1.2) has two primary singularities at the origin blowing expansion, where $P_1$ is a saddle point and $P_2$ is an unstable node. When $s > m$, the two primary singularities blown out by the system are both saddle points, so it is known that the origin of the system (1.2) is an unstable point.

## 2.2 Equilibrium point analysis

Discuss the equilibrium point on $\Omega = \{(x, y) | x > 0, y \geq 0\}$ for solving (1.2), i.e., solving the following system of equations:

$$\begin{cases} \dfrac{x}{1+\lambda y}(1-x)(x-m) - axy = 0 \\ sy\left(1-\dfrac{y}{x}\right) = 0 \end{cases} \quad (2.2)$$

Solving the system of equations (2.2) leads to the following conclusion:

**Theorem 2:** The system has two boundary equilibria $E_1 = (1,0)$ and $E_2 = (m,0)$.

**Proof:** from (2.2), the boundary equilibrium point of the system occurs when $y = 0$. Therefore, when $y = 0$, it is obtained that:

$$(1-x)(x-m) = 0 \quad (2.3)$$

That is, we can get $x_1 = 1$, $x_2 = m$.

**Theorem 3:** Noting $a_1^* = m+1-2\sqrt{m}$, $\Delta = (m+1-a)^2 - 4m(1+\lambda a)$, the internal equilibrium point of the system (1.2) results in the following:

(1) The system (1.2) has no internal equilibrium point when one of the following conditions is satisfied:

(i). $0 < m < 1$, $a_1^* \leq a < m+1$.

(ii). $0 < m < 1$, $0 < a < a_1^*$, $\lambda > \dfrac{a^2 - 2(m+1)a + (m-1)^2}{4ma}$.

(2) When $0 < a < a_1^*$, $\lambda = \dfrac{a^2 - 2(m+1)a + (m-1)^2}{4ma}$, the system has a unique internal equilibrium point

$$E_3 = (x_3, y_3) = \left(\dfrac{m+1-a}{2(1+\lambda a)}, \dfrac{m+1-a}{2(1+\lambda a)}\right) = \left(\dfrac{2m}{1+m-a}, \dfrac{2m}{1+m-a}\right).$$

(3) When $0 < m < 1$, $0 < a < a_1^*$, $0 < \lambda < \dfrac{a^2 - 2(m+1)a + (m-1)^2}{4ma}$, The system has two internal equilibrium points, which are

$$E_4 = (x_4, y_4) = \left(\dfrac{1+m-a-\sqrt{\Delta}}{2(1+\lambda a)}, \dfrac{1+m-a-\sqrt{\Delta}}{2(1+\lambda a)}\right),$$

$$E_5 = (x_5, y_5) = \left( \frac{1+m-a+\sqrt{\Delta}}{2(1+\lambda a)}, \frac{1+m-a+\sqrt{\Delta}}{2(1+\lambda a)} \right).$$

**Proof:** Solve for the positive equilibrium point i.e. $x > 0$ and $y > 0$. So from the second equation of (2.2) we have:

$$y = x \tag{2.4}$$

Bringing equation (2.4) into the first equation of (2.2) yields:

$$(1+\lambda a)x^2 - (1+m-a)x + m = 0 \tag{2.5}$$

Since $1+\lambda a > 0$ and $m > 0$, by the quadratic property, the equation (2.5) has no positive roots when $1+m-a \leq 0$. So only the case when $1+m-a > 0$, i.e. $0 < a < 1+m$, needs to be discussed.

Denote the discriminant of equation (2.5) as:

$$\Delta = (m+1-a)^2 - 4m(1+\lambda a)$$

Since $m \in (0,1)$, consider $\Delta$ subject to the range of values of $m$:

(a) Let $\Delta < 0$, we can get $\lambda > \dfrac{a^2 - 2(m+1)a + (m-1)^2}{4ma}$, because $ma > 0$, consider the sign of $f(a) = a^2 - 2(m+1)a + (m-1)^2$, and write $\Delta_a = 4(1+m)^2 - 4(m-1)^2$.

For $\forall m \in (0,1)$, we have $f(0) = (m-1)^2 > 0$, the axis of symmetry $x' = 2(m+1) > 0$ of $f(a)$, and $\Delta_a = 16m > 0$, so that two positive roots are obtained when $f(a) = 0$:

$$a_i^* = m+1 \pm \sqrt{m}, i = 1, 2$$

and with $a_1^* < m+1 < a_2^*$.

So when $0 < a < a_1^*$, $f(a) > 0$; when $a = a_1^*$, $f(a) = 0$; when $a_1^* < a < m+1$, $f(a) < 0$.

According to the sign of $f(a)$, there is $\Delta < 0$ and equation (2.5) has no positive roots when one of the following conditions is satisfied, i.e., the system (1.2) has no internal equilibrium point:

(i). $0 < m < 1$, $a_1^* \leq a < m+1$.

(ii). $0 < m < 1$, $0 < a < a_1^*$, $\lambda > \dfrac{a^2 - 2(m+1)a + (m-1)^2}{4ma}$.

(b) Let $\Delta = 0$, From (a), it follows that when $0 < m < 1$, $0 < a < a_1^*$,

$$\lambda = \frac{a^2 - 2(m+1)a + (m-1)^2}{4ma}, \Delta = 0,$$ The equation (2.5) has a unique positive root $x_3$, and

$$x_3 = \frac{m+1-a}{2(1+\lambda a)} = \frac{2m}{1+m-a},$$ So the system (1.2) has a unique internal equilibrium point

$E_3 = (x_3, y_3)$, where $y_3 = x_3$.

(c) Let $\Delta > 0$, From (a), it follows that when $0 < m < 1$, $0 < a < a_1^*$,

$$0 < \lambda < \frac{a^2 - 2(m+1)a + (m-1)^2}{4ma}$$ ,Equation (2.5) has two positive roots

$$x_4 = \frac{1+m-a-\sqrt{\Delta}}{2(1+\lambda a)} \text{ and } x_5 = \frac{1+m-a+\sqrt{\Delta}}{2(1+\lambda a)},$$ and $x_4 < x_5$, The system (1.2) has two internal

equilibria $E_4 = (x_4, y_4)$ and $E_5 = (x_5, y_5)$, where $y_i = x_i, i = 4, 5$.

## 2.3 Equilibrium point dynamical form

The dynamic form at the equilibrium point is the key to further study of the model, and the analysis of the stability of the system at the equilibrium point allows to propose solutions to the practical problems corresponding to the model. For this purpose, the eigenvalues of the Jacobi matrix of the system at the equilibrium point need to be considered. If $E = (x, y)$ is the equilibrium point of the system (1.2), then the Jacobi matrix of this system at the equilibrium point is:

$$J(E) = \begin{pmatrix} \frac{(1-2x)(x-m) + x(1-x)}{1+\lambda y} - ay & -\frac{\lambda x(1-x)(x-m)}{(1+\lambda y)^2} - ax \\ s & -s \end{pmatrix}$$

Based on the Jacobi matrix at the equilibrium point, the analysis of its eigenvalues leads to the following theorem:

**Theorem 4:** The boundary equilibrium point $E_1 = (1, 0)$ of the system (1.2) is the saddle point.

**Proof:** the Jacobian matrix of the system (1.2) at the equilibrium point $E_1$ is:

$$J(E_1) = \begin{pmatrix} m-1 & -a \\ 0 & s \end{pmatrix}$$

The eigenvalues of $J(E_1)$ are $\lambda_1 = m-1$, $\lambda_2 = s$. For the strong Allee effect, we have $0 < m < 1$, so we have $\lambda_1 < 0, \lambda_2 > 0$. Therefore, the boundary equilibrium point $E_1$ is the saddle point.

**Theorem 5:** The boundary equilibrium point $E_2 = (m, 0)$ of the system (1.2) is an unstable node.

**Proof:** the Jacobian matrix of the system (1.2) at the equilibrium point $E_2$ is:

$$J(E_2) = \begin{pmatrix} m(1-m) & -am \\ 0 & s \end{pmatrix}$$

The eigenvalues of $J(E_2)$ are $\lambda_1 = m(1-m)$, $\lambda_2 = s$. For the strong Allee effect, we have $0 < m < 1$, so we have $\lambda_1 > 0, \lambda_2 > 0$. Therefore, the boundary equilibrium point $E_2$ is an unstable node.

**Theorem 6:** When $0 < m < 1$, $0 < a < a_1^*$, $\lambda = \dfrac{a^2 - 2(m+1)a + (m-1)^2}{4ma}$, for the internal equilibrium point $E_3 = (x_3, y_3)$ of the system (1.2), denoted $s_0 = \dfrac{-2x_3^2 + (m+1)x_3}{1 + \lambda y_3}$, we have the following conclusion:

(1) When $0 < s < s_0$, $E_3$ is the excluded saddle node.

(2) When $s > s_0$, $E_3$ is the attracted saddle node

(3) When $s = s_0$, $E_3$ is the cusp with residual dimension 2.

**Proof:** the Jacobian matrix of the system (1.2) at the internal equilibrium point $E_3$ is:

$$J(E_3) = \begin{pmatrix} \dfrac{(1-2x_3)(x_3-m) + x_3(1-x_3)}{1+\lambda y_3} - ay_3 & -\dfrac{\lambda x_3(1-x_3)(x_3-m)}{(1+\lambda y_3)^2} - ax_3 \\ s & -s \end{pmatrix}$$

Since the equilibrium point $E_3$ satisfies equation (2.5), i.e:

$$(1 + \lambda a)x_3^2 = (1 + m - a)x_3 - m \tag{2.6}$$

Let $tr(J(E_3))$ denote the trace of matrix $J(E_3)$ and $det(J(E_3))$ denote the determinant of matrix $J(E_3)$. Then:

$$tr(J(E_3)) = \dfrac{-2x_3^2 + (m+1)x_3}{1 + \lambda y_3} - s = s_0 - s$$

and

$$det(J(E_3)) = -s\left(\frac{(1-2x_3)(x_3-m)+x_3(1-x_3)}{1+\lambda y_3} - ay_3 - \frac{\lambda x_3(1-x_3)(x_3-m)}{(1+\lambda y_3)^2} - ax_3\right) = \frac{s\left((1+\lambda a)x_3^2 - m\right)}{1+\lambda y_3}$$

According to $\Delta = 0$, we can obtain $2m = \dfrac{(1+m-a)^2}{2(1+\lambda a)}$, and combining with equation (2.6), we can obtain:

$$det(J(E_3)) = \frac{s\left((1+\lambda a)x_3^2 - 2m\right)}{1+\lambda y_3} = 0$$

Next, consider the notation of $s_0$, from $0 < m < 1, 0 < a < a_1^*$, $x_3 = \dfrac{2m}{1+m-a}$, we can obtain:

$$(m+1) - 2x_3 = m + 1 - \frac{4m}{1+m-a} > m + 1 - 2\sqrt{m} = \left(\sqrt{m}-1\right)^2 > 0 \quad (2.7)$$

So from equation (2.7) it follows that $s_0 = \dfrac{-2x_3^2 + (m+1)x_3}{1+\lambda y_3} > 0$.

(1) When $s \neq s_0$, $tr(J(E_3)) \neq 0$, The matrix $J(E_3)$ has a zero eigenvalue. So doing the transformation:

$$\begin{cases} u = x - x_3 \\ v = y - y_3 \end{cases}$$

Translating the boundary point to the origin and performing a Taylor expansion at the origin yields:

$$\begin{cases} \dot{u} = a_1 u + a_2 v + a_3 u^2 + a_4 uv + a_5 v^2 + O\left((u,v)^2\right) \\ \dot{v} = b_1 u + b_2 v + b_3 u^2 + b_4 uv + b_5 v^2 + O\left((u,v)^2\right) \end{cases} \quad (2.8)$$

Where:

$$a_1 = \frac{-2x_3^2 + (m+1)x_3}{1+\lambda y_3}, \quad a_2 = -\frac{\lambda x_3(1-x_3)(x_3-m)}{(1+\lambda y_3)^2} - ax_3$$

$$a_3 = \frac{-3x_3 + (m+1)}{1+\lambda y_3}, \quad a_4 = -\frac{\lambda\left((1-2x_3)(x_3-m)+x_3(1-x_3)\right)}{(1+\lambda y_3)^2} - a$$

$$a_5 = \frac{2a\lambda^2 x_3^2}{(1+\lambda y_3)^2}, \quad b_1 = s, \quad b_2 = -s, \quad b_3 = -\frac{s}{x_3}, \quad b_4 = \frac{2s}{x_3}, \quad b_5 = -\frac{s}{x_3}。$$

From $det(J(E_3)) = a_1 b_2 - a_2 b_1 = 0$, we can get $a_1 = -a_2$.

To transform the system into a standard type, do the transformation:

$$\begin{pmatrix} u \\ v \end{pmatrix} = \begin{pmatrix} a_2 & a_1 \\ -a_1 & b_1 \end{pmatrix} \begin{pmatrix} \tilde{x} \\ \tilde{y} \end{pmatrix}$$

and doing a time transformation $dt = (s_0 - s)d\tau$, the system (2.8) can be written as:

$$\begin{cases} \dot{\tilde{x}} = c_1 \tilde{x}^2 + c_2 \tilde{x}\tilde{y} + c_3 \tilde{y}^2 + O\left((\tilde{x}, \tilde{y})^2\right) \\ \dot{\tilde{y}} = d_1 \tilde{y} + d_2 \tilde{x}^2 + d_3 \tilde{x}\tilde{y} + d_4 \tilde{y}^2 + O\left((\tilde{x}, \tilde{y})^2\right) \end{cases} \quad (2.9)$$

Where:

$$c_1 = \frac{s_0 - s}{a_2}\left(a_2^2 a_3 - a_1 a_2 a_4 + a_1^2 a_5\right) + a_2^2(a_3 - b_3) - a_1 a_2(a_4 - b_4) + a_1^2(a_5 - b_5)$$

$$c_2 = \frac{s_0 - s}{a_2}\left(2a_1 a_2 a_3 + a_4(a_2 b_1 - a_1^2) - 2a_1 b_1 a_5\right) + 2a_1 a_2(a_3 - b_3) + (a_2 b_1 - a_1^2)(a_4 - b_4) - 2a_1 b_1(a_5 - b_5)$$

$$c_3 = \frac{s_0 - s}{a_2}\left(a_1^2 a_3 + a_1 a_4 b_1 + a_5 b_1^2\right) + a_1^2(a_3 - b_3) + a_1 b_1(a_4 - b_4) - b_1^2(a_5 - b_5)$$

$$d_1 = (a_1 - b_1)^2, \quad d_2 = a_2^2(a_3 - b_4) - a_1 a_2(a_4 - b_4)$$

$$d_3 = 2a_1 a_2(a_3 - b_3) + (a_2 b_1 - a_1^2)(a_4 - b_4) - 2a_1 b_1(a_5 - b_5)$$

$$d_4 = a_1^2(a_3 - b_3) + a_1 b_1(a_4 - b_4) - b_1^2(a_5 - b_5)$$

From equation (2.9), we get that the coefficient of $\tilde{x}^2$ is:

$$c_1 = -a_1(s_0 - s)(a_3 + a_4 + a_5) + a_1^2(a_3 - b_3 + a_4 - b_4 + a_5 - b_5) \neq 0$$

Therefore, from the literature [9], the equilibrium point $E_3$ is a saddle node at this time. And when $0 < s < s_0$, $c_1 < 0$, and the time transformation is positive, so the equilibrium point $E_3$ is a repulsive saddle node at this time. When $s > s_0$, $c_1 > 0$, and the time transformation is negative, so the equilibrium point $E_3$ is the saddle node of the attraction.

(2) When $s = s_0$, $tr(J(E_3)) = det(J(E_3)) = 0$, The eigenvalues of the matrix $J(E_3)$ are all zero eigenvalues. Therefore, the system (2.8) is reduced to the standard type by doing the transformation:

$$\begin{pmatrix} u \\ v \end{pmatrix} = \begin{pmatrix} a_2 & 0 \\ -a_1 & 1 \end{pmatrix} \begin{pmatrix} X \\ Y \end{pmatrix}$$

Then system (2.8) can be written as:

$$\begin{cases} \dot{X} = Y + e_1 X^2 + e_2 XY + e_3 Y^2 + O\left((X,Y)^3\right) \\ \dot{Y} = f_1 X^2 + f_2 XY + f_3 Y^2 + O\left((X,Y)^3\right) \end{cases} \quad (2.10)$$

Where:

$$e_1 = a_2 a_3 - a_1 a_4 + \frac{a_1^2 a_5}{a_2}, \quad e_2 = a_4 - \frac{2 a_1 a_5}{a_2}, \quad e_3 = \frac{a_5}{a_2}$$

$$f_1 = -a_2^2 (a_3 - b_3) + a_1 a_2 (a_4 - b_4) - a_1^2 (a_5 - b_5)$$

$$f_2 = -a_2 (a_4 - b_4) - 2 a_1 (a_5 - b_5), \quad f_3 = -(a_5 - b_5)$$

From the literature [10], the system (2.10) is equivalent to the system (2.11) in the vicinity of the origin, i.e:

$$\begin{cases} \dot{X} = Y + O\left((X,Y)^2\right) \\ \dot{Y} = f_1 X^2 + (f_2 + 2e_1) XY + O\left((X,Y)^3\right) \end{cases} \quad (2.11)$$

Since:

$$f_1 = -a_2^2 (a_3 - b_3) + a_1 a_2 (a_4 - b_4) - a_1^2 (a_5 - b_5) = -b_1^2 (a_3 + a_4 + a_5) \neq 0$$

$$f_2 + 2e_1 = -b_1 (2 a_3 + a_4)$$

From Theorem 3 in the literature [11], the equilibrium point $E_3$ is a cusp with a residual dimension of 2 when $f_1 (f_2 + 2e_1) \neq 0$.

**Theorem 7:** When $0 < m < 1$, $0 < a < a_1^*$, $\lambda = \dfrac{a^2 - 2(m+1)a + (m-1)^2}{4ma}$, for the internal equilibrium point $E_4 = (x_4, y_4)$ of the system (1.2), denoted $s^* = \dfrac{-2 x_4^2 + (m+1) x_4}{1 + \lambda y_4}$, we have the following conclusion:

(1) When $0 < s < s^*$, $E_4$ E is an unstable node or focus.

(2) 当 $s > s^*$, $E_4$ is a stable node or focus.

(3) 当 $s = s^*$ 时, $E_4$ is the weak center.

**Proof:** the Jacobian matrix of the system (1.2) at the internal equilibrium point $E_4$ is:

$$J(E_4) = \begin{pmatrix} \dfrac{(1-2x_4)(x_4-m)+x_4(1-x_4)}{1+\lambda y_4} - ay_4 & -\dfrac{\lambda x_4(1-x_4)(x_4-m)}{(1+\lambda y_4)^2} - ax_4 \\ s & -s \end{pmatrix}$$

Since the equilibrium point E satisfies equation (2.5), i.e:

$$(1+\lambda a)x_4^2 = (1+m-a)x_4 - m \tag{2.12}$$

Also according to $\Delta > 0$, then $0 < \lambda < \dfrac{a^2 - 2(m+1)a + (m-1)^2}{4ma}$, combined with equation (2.12), we have:

$$det(J(E_4)) = s\left( \dfrac{(x_4-m)-x_4(1-x_4)}{1+\lambda y_4} + \dfrac{\lambda x_4(1-x_4)(x_4-m)}{(1+\lambda y_4)^2} \right) > 0$$

and

$$tr(J(E_4)) = \dfrac{-2x_4^2 + (m+1)x_4}{1+\lambda y_4} - s = s^* - s$$

Next consider the sign of $s^*$. Since:

$$(m+1) - 2x_4 = m+1 - \dfrac{1+m-a-\sqrt{\Delta}}{1+\lambda a} > a + \sqrt{\Delta} > 0 \tag{2.13}$$

So from equation (2.13) we know that, $s^* = \dfrac{-2x_4^2 + (m+1)x_4}{1+\lambda y_4} > 0$.

Therefore, when $0 < s < s^*$, there are $tr(J(E_4)) > 0$, $det(J(E_4)) > 0$. Then $E_4$ is an unstable focus or node.

When $s > s^*$, there are $tr(J(E_4)) < 0$, $det(J(E_4)) > 0$. Then $E_4$ is a stable focus or node.

When $s = s^*$, there are $tr(J(E_4)) = 0$. The eigenvalues are $\lambda_{1,2} = \pm wi$, where

$$w = \dfrac{\sqrt{4det(J(E_4)) - (tr(J(E_4)))^2}}{2},$$ and at this point $E_4$ is the weak center.

**Theorem 8:** When $0 < m < 1$, $0 < a < a_1^*$, $\lambda = \dfrac{a^2 - 2(m+1)a + (m-1)^2}{4ma}$, The internal equilibrium point $E_5 = (x_5, y_5)$ of the system (1.2) is a stable node or focal point.

**Proof:** the Jacobian matrix of the system (1.2) at the internal equilibrium point $E_5$ is:

$$J(E_5) = \begin{pmatrix} \dfrac{(1-2x_5)(x_5-m)+x_5(1-x_5)}{1+\lambda y_5} - ay_5 & -\dfrac{\lambda x_5(1-x_5)(x_5-m)}{(1+\lambda y_5)^2} - ax_5 \\ s & -s \end{pmatrix}$$

From Theorem 7, it follows that:

$$det(J(E_5)) = s\left(\dfrac{(x_5-m)-x_5(1-x_5)}{1+\lambda y_5} + \dfrac{\lambda x_5(1-x_5)(x_5-m)}{(1+\lambda y_5)^2}\right) > 0$$

and

$$tr(J(E_5)) = \dfrac{-2x_5^2+(m+1)x_5}{1+\lambda y_5} - s$$

From $0 < \lambda < \dfrac{a^2-2(m+1)a+(m-1)^2}{4ma}$ and $0 < m < 1$ we have:

$$(m+1) - 2x_5 < \dfrac{(m-1)(m+1-a)^2 - 4m\sqrt{\Delta}}{(m+1-a)^2} < 0 \quad (2.14)$$

So from equation (2.14) we know that $\dfrac{-2x_5^2+(m+1)x_5}{1+\lambda y_5} < 0$, i.e. $tr(J(E_5)) < 0$. So $E_5$ is a stable node or focus.

## 3  Bifurcation of the Leslie-Gower predatory-prey model

The study of the dynamical behaviors such as bifurcation of nonlinear predator-prey systems is beneficial to the rational conservation of natural resources by humans. This subsection analyzes the possible bifurcations of the system (1.2), including saddle-node bifurcation and Hopf bifurcation, and gives the conditions for the existence of these bifurcations based on the results of the dynamic shape analysis of the equilibrium point in Section 2.

### 3.1  Saddle-node bifurcation

It follows from Theorem 2 that when $0 < m < 1$, $0 < a < a_1^*$, if

$0 < \lambda < \dfrac{a^2-2(m+1)a+(m-1)^2}{4ma}$, the system has two internal equilibria $E_4$ and $E_5$; If

$\lambda > \dfrac{a^2-2(m+1)a+(m-1)^2}{4ma}$, the system has no internal equilibrium point; if

$\lambda = \dfrac{a^2-2(m+1)a+(m-1)^2}{4ma}$, there is a unique positive equilibrium point $E_3$.

So let $\lambda = \lambda_{SN} = \dfrac{a^2 - 2(m+1)a + (m-1)^2}{4ma}$, when $\lambda$ varies around $\lambda_{SN}$, the number of boundary equilibrium points changes, giving rise to saddle-node branches.

**Theorem 9:** When $s \neq s_0$, the system (1.2) undergoes a saddle-node bifurcation at $E_3$ with a critical bifurcation parameter $\lambda = \lambda_{SN} = \dfrac{a^2 - 2(m+1)a + (m-1)^2}{4ma}$.

**Proof:** By Theorem 6, the Jacobian matrix of the system (1.2) at $E_3 = (x_3, y_3)$ is:

$$J(E_3) = \begin{pmatrix} \dfrac{(1-2x_3)(x_3-m)+x_3(1-x_3)}{1+\lambda y_3} - ay_3 & -\dfrac{\lambda x_3(1-x_3)(x_3-m)}{(1+\lambda y_3)^2} - ax_3 \\ s & -s \end{pmatrix}$$

Because $s \neq s_0$, so $J(E_3)$ has a zero eigenvalue, let $v$, $w$ be the eigenvectors corresponding to the zero eigenvalues of $J(E_3)$ and $J^T(E_3)$, respectively. Then we can obtain:

$$v = \begin{pmatrix} v_1 \\ v_2 \end{pmatrix} = \begin{pmatrix} 1 \\ 1 \end{pmatrix}$$

$$w = \begin{pmatrix} w_1 \\ w_2 \end{pmatrix} = \begin{pmatrix} s \\ -\dfrac{\lambda a x_3}{1+\lambda y_3} - ax_3 \end{pmatrix}$$

Let

$$F(x, y) = \begin{pmatrix} \dot{x} \\ \dot{y} \end{pmatrix} = \begin{pmatrix} F_1 \\ F_2 \end{pmatrix} = \begin{pmatrix} \dfrac{x}{1+\lambda y}(1-x)(x-m) - axy \\ sy\left(1-\dfrac{y}{x}\right) \end{pmatrix}$$

We obtain:

$$F_a(E_3, \lambda_{SN}) = \begin{pmatrix} -x_3 y_3 \\ 0 \end{pmatrix}$$

$$D^2 F(E_3, \lambda_{SN})(v,v) = \begin{pmatrix} \dfrac{\partial^2 F_1}{\partial x^2} v_1^2 + 2\dfrac{\partial^2 F_1}{\partial x \partial y} v_1 v_2 + \dfrac{\partial^2 F_1}{\partial y^2} v_2^2 \\ \dfrac{\partial^2 F_2}{\partial x^2} v_1^2 + 2\dfrac{\partial^2 F_2}{\partial x \partial y} v_1 v_2 + \dfrac{\partial^2 F_2}{\partial y^2} v_2^2 \end{pmatrix} = \begin{pmatrix} -\dfrac{2\lambda_{SN}^2 a x_3^2}{(1+\lambda_{SN} y_3)^2} \end{pmatrix}$$

So we get:

$$w^T F_a(E_3, \lambda_{SN}) = -sx_3 y_3 \neq 0$$

$$w^T\left[D^2F(E_3,\lambda_{SN})(v,v)\right]=-\frac{2\lambda_{SN}^2 ax_3^2}{(1+\lambda_{SN}y_3)^2}\neq 0$$

By Sotomayor's theorem, it is known that the system (1.2) has a saddle-node bifurcation at the equilibrium point $E_3$.

### 3.2 Hopf bifurcation

By Theorem 7, when $0<m<1$, $0<a<a_1^*$, $0<\lambda<\dfrac{a^2-2(m+1)a+(m-1)^2}{4ma}$, if $0<s<s^*$, then $E_4$ is an unstable focus or node; if $s>s^*$, then $E_4$ is a stable focus or node; if $s=s^*$, then $E_4$ is the weak center. So observe whether the system generates Hopf bifurcation when the parameter $s$ is changed.

**Theorem 10:** When $s=s^*$, the system (1.2) has a Hopf bifurcation at the equilibrium point $E_4$.

**Proof:** Let

$$\mu(s)=tr(J(E_4))=\frac{-2x_4^2+(m+1)x_4}{1+\lambda y_4}-s \tag{3.1}$$

So we can get $\mu(s^*)=0$ and $\mu'(s^*)<0$, so the system will have a Hopf bifurcation at $E_4$.

Next, the first Lyapunov coefficient is calculated and the direction of the bifurcation is determined. When $s=s^*$, the equilibrium point is translated to the origin and the coordinate transformation is chosen:

$$\begin{cases} x=x_4+\xi_1 \\ y=y_4+\xi_2 \end{cases}$$

and do Taylor expansion at the origin to obtain the system:

$$\begin{cases} \dot{\xi}_1 = m_1\xi_1+m_2\xi_2+m_3\xi_1^2+m_4\xi_1\xi_2+m_5\xi_2^2+m_6\xi_1^3+m_7\xi_1^2\xi_2+m_8\xi_1\xi_2^2+m_9\xi_2^3+O\left((\xi_1,\xi_2)^3\right) \\ \dot{\xi}_2 = n_1\xi_1+n_2\xi_2+n_3\xi_1^2+n_4\xi_1\xi_2+n_5\xi_2^2+n_6\xi_1^3+n_7\xi_1^2\xi_2+n_8\xi_1\xi_2^2+O\left((\xi_1,\xi_2)^3\right) \end{cases}$$

Where:

$$m_1=\frac{(1-2x_4)(x_4-m)+x_4(1-x_4)}{1+\lambda y_4}-ay_4, \quad m_2=-\frac{\lambda x_4(1-x_4)(x_4-m)}{(1+\lambda y_4)^2}-ax_4$$

$$m_3=\frac{-(x_4-m)+x_4(1-x_4)}{1+\lambda y_4}, \quad m_4=-\frac{\lambda\left((1-2x_4)(x_4-m)+x_4(1-x_4)\right)}{(1+\lambda y_4)^2}-a$$

$$m_5=\frac{\lambda^2 x_4(1-x_4)(x_4-m)}{(1+\lambda y_4)^3}-ax_4, \quad m_6=-\frac{1}{1+\lambda y_4}, \quad m_7=\frac{\lambda\left((x_4-m)-(1-2x_4)\right)}{(1+\lambda y_4)^3},$$

$$m_8 = \frac{\lambda^2((1-2x_4)(x_4-m)+x_4(1-x_4))}{(1+\lambda y_5)^3}, \quad m_9 = -\frac{\lambda^3 x_4(1-x_4)(x_4-m)}{(1+\lambda y_4)^4}$$

$$n_1 = s^*, \quad n_2 = -s^*, \quad n_3 = -\frac{s^*}{x_4^2}, \quad n_4 = \frac{2s^*}{x_4}, \quad n_5 = -\frac{s^*}{x_4}, \quad n_6 = \frac{s^*}{x_4^3}, \quad n_7 = -\frac{s^*}{x_4^3}, \quad n_8 = \frac{2s^*}{x_4^2}.$$

From $tr(J(E_4)) = m_1 + n_2 = 0$, we can get $m_1 = -n_2$.

According to the theorem in Perko [11], the first Lyapunoov coefficient can be obtained as:

$$l_1 = \frac{-3\pi}{2m_2(m_1n_2 - m_2n_1)^{3/2}}\{m_1n_1\psi_1 + m_1m_2\psi_2 + n_1^2\psi_3 - 2m_1n_1\psi_4 - 2m_1m_2\psi_5 - m_2^2\psi_6 - (m_2n_1 - 2m_1^2)\psi_7 - (m_1^2 + m_2n_1)\psi_8\}$$

Where:

$$\psi_1 = m_4^2 + m_4n_5 + m_5n_4, \quad \psi_2 = n_4^2 + m_3n_4 + m_4n, \quad \psi_3 = m_4m_5 + 2m_5n_5$$

$$\psi_4 = n_5^2 - m_3m_5, \quad \psi_5 = m_3^2 - n_3n_5, \quad \psi_6 = 2m_3n_3 + n_3n_4, \quad \psi_7 = n_4n_5 - m_4m_3$$

$$\psi_8 = -3m_6m_2 + 2m_1(m_7 + n_8) + (n_1n_8 - m_2n_7)$$

Therefore, according to the expression of $l_1$, the following theorem can be obtained:

**Theorem 11:** For system (1.2)

(1) When $l_1 < 0$, The system (1.2) has supercritical Hopf bifurcation generation near $E_4$.

(2) When $l_1 > 0$, The system (1.2) has a subcritical Hopf bifurcation generated near $E_4$.

# 4 Conclusion

In this paper, a Leslie-Gower predatory-prey model with fear effect and Allee effect is studied, and the origin is first obtained as an unstable point using the blow-up method. Subsequently, by the analysis of the equilibrium points, the model (1.2) always has two boundary equilibrium points, which are the saddle point and the unstable node. And the existence of one or two internal positive equilibrium points of the system under certain conditions. Through theoretical analysis, the fear effect and Allee effect have a great influence on the dynamic properties of the model, and the stable positive equilibrium point may become unstable with the change of the fear effect and Allee effect. When the intrinsic growth rate of the predator reaches a critical value $s^*$, the system generates a Hopf bifurcation at the internal equilibrium point.

[1]College of Mathematics and Statistics, Chongqing University, Chongqing, 401331, P. R. China
Email address: chenlinchn@126.com
[*]College of Mathematics and Statistics, Chongqing University, Chongqing, 401331, P. R. China
Email address: zhuchangrong126@126.com